\documentclass[reqno, 12pt]{amsart}
\usepackage[all]{xy}
 \RequirePackage{amsgen}
\RequirePackage{amstext} \RequirePackage{amsbsy}
\RequirePackage{amsopn} \RequirePackage{amsthm}
\RequirePackage{amssymb} \RequirePackage{epsfig}
\RequirePackage{amscd}
\usepackage{amsmath, latexsym, amssymb, fancyhdr}

\usepackage{graphics}
\usepackage{wrapfig}

\usepackage[mathscr]{eucal}

\newtheorem{lemma}{Lemma}

\newtheorem{prop}{Proposition}

\newtheorem{theorem}{Theorem}

\theoremstyle{definition}
\newtheorem{defn}{Definition}

\theoremstyle{remark}
\newtheorem{remark}{Remark}

\newtheorem{example}{Example}

\let\a\alpha    \let\c\gamma  
  \let\g\gamma
  \let\l\lambda   
      
\let\GL\Lambda
\def\C{\mathbb C}

\def\E{\mathcal E}

\def\fH{\frak H}

\def\GL{\rm{GL}}

\def\Tr{{\rm Tr}}
\def\Gal{{\rm Gal}}

\def\Z{{\mathbb Z}}

\def\Q{{\mathbb Q}}

\def\G{\Gamma}
\def\ord{{\rm ord}}

\def\SL{SL_2(\mathbb Z)}

\title[\textit{Modular subgroups and modular forms}]{Finite index subgroups of the modular group  and their modular forms}

\author{Ling Long}

\address{Department of Mathematics\\Iowa State University\\Ames, IA 50011 \\USA}
\email{linglong@iastate.edu}

\thanks{ The author was
supported in part by an NSF-AWM mentoring travel grant for women.
She would like to thank the Pennsylvania State University  for
hosting her visit during May-June 2006.}

\begin{document}

\begin{abstract}
   Classically, congruence  subgroups of the modular group, which can be described by congruence relations, play important roles
  in group theory and modular forms. In reality, the majority of finite index  subgroups of the modular group are noncongruence. These groups as well as their modular forms
  are central players of this survey article.  Differences
  between congruence and noncongruence  subgroups and modular forms  will be
  discussed.  We will mainly focus on three interesting
  aspects of modular forms for noncongruence  subgroups:
  the unbounded denominator property,  modularity of the Galois representation arising from
  noncongruence cuspforms, and Atkin and Swinnerton-Dyer
  congruences.
\end{abstract}

\maketitle
\section{Introduction}The special linear group $SL_2(\Z)$ consists
of all 2-by-2 integral matrices with determinant 1. It is one of the
most well-known and important discrete groups.   A finite index
subgroup  of the modular group is said to be \emph{congruence }if it
contains the kernel of the natural modulo $N$ homomorphism from
$SL_2(\Z)$ to $SL_2(\Z/N\Z)$ for some positive integer $N$;
otherwise, it is called a \emph{noncongruence} subgroup. The existence of
noncongruence subgroups of $SL_2(\Z)$ was confirmed by Fricke
\cite{fricke1886} and Pick \cite{Pick1886}. More noncongruence
subgroups were constructed by Reiner \cite{reiner58}, Newman
\cite{newman65}, Rankin \cite{rankin67}, etc. In contrast, any
finite index subgroup of $SL_n(\Z)$ with $n\ge 3$ is congruence
\cite{bass-lazard-serre64}. A famous theorem of Bely\u{\i} implies
any compact complex smooth irreducible curve defined over
$\overline{\Q}$, the algebraic closure of $\Q$, can be realized (in
many ways) as a modular curve (defined in \ref{subsec:modcurve}) for
a finite index subgroup of $\SL$ \cite{belyi79}. Consequently,
finite index subgroups arise naturally in many fields such as the
theory of dessin d'enfant which includes Galois coverings of the
projective line \cite{birchb94}, triangular groups
\cite{jon79,Jones-Sigerman94}, as well as  the theory of translation
surfaces \cite{H-L05}. As a matter of fact, noncongruence subgroups
of $\SL$ predominate congruence subgroups (Section 2.3 below).

Modular forms are generalizations of elliptic functions (doubly
periodic functions) and enjoy plenty of wonderful properties. Many
important one variable formal power series  arising in combinatorics
or physics, such as the partition function, turn out to be closely
related to modular forms.

The theory of congruence (elliptic) modular forms is now
well-developed due to the efforts of many mathematicians. It has
very broad impacts including some important applications to
information technology \cite{sarnak90,li1}.  In the early days,
computational observations and conjectures were significant driving
forces for the development of this theory. For example, Ramanujan
observed coefficients of $q\prod_{n\ge 1}(1-q^n)^{24}=\sum \tau(n)
q^n$ satisfy special three-term recursions for small primes $p$
\begin{equation}\label{eq:Delta}
\tau(np)-\tau(p)\tau(n)+p^{11}\tau(n/p)=0, \forall n\in \mathbb N,
\end{equation}where $\tau(n/p)=0$ if $p\nmid n$.
In  \cite{mordell17} Mordell showed that these recursions hold for
all primes $p$ by using  essentially Hecke operators. In
\cite{hecke37} Hecke defined  Hecke operators
 and proved that  they satisfy
$(2d+1)$-term relations on each $d$-dimensional space of congruence
modular forms. Shortly after,
  Petersson introduced the Petersson inner product
 via which he   showed  the Hecke operators can be diagonalized
simultaneously \cite{petersson391}. The readers are referred to
\cite[Section 2]{lly052} for a summary of the theory of newforms
(congruence modular form with particular properties) developed by
Atkin and Lehner \cite{a-l70}, Miyake \cite{miy71}, and Li
\cite{wli75}. The current considerable interest in modular forms is
largely due to the vital role  they played in the celebrated proof
of  Fermat's Last Theorem (FLT) by Wiles \cite{wil95}. In his proof,
Wiles used a result of Frey and Ribet which says  the
Taniyama-Shimura-Weil (TSW) conjecture implies FLT. The TSW
conjecture asserts the $l$-adic Galois representations attached to
an elliptic curve (a genus one complete curve) defined over $\Q$ are
isomorphic to the $l$-adic Deligne representations attached to a
certain newform. The final proof of TSW (\cite{wil95,tw95,bcdt}) is
a major breakthrough in the very influential Langlands program which
is a system of conjectures connecting the
 $l$-adic cohomology  theory of algebraic varieties  with representation theory of cuspidal automorphic
 forms \cite{Gelbart84, Murty-langlands}. In the sequel, we say an $l$-adic Galois
representation is ``modular" when it  arises from a congruence
automorphic form.
\smallskip

Modular forms for noncongruence subgroups have already been
considered by Fricke and Klein. However, there was no serious
investigation of these functions till the late 1960's as it was
commonly believed that the Hecke theory was missing. Indeed, Atkin
conjectured that Hecke operators defined classically have trivial
actions on noncongruence modular forms and his conjecture was proven
in part by Serre \cite{tho89} and Scholl \cite{sch97} and in general
by Berger \cite{berger94}.

However, the investigation in the past 40 years  reveals  the
mathematical structure of noncongruence modular forms is still
wonderfully rich. In \cite{a-sd}, Atkin and Swinnerton-Dyer  used
computers together with algebra and combinatorics to work with
noncongruence subgroups  and noncongruence modular forms. They
brought up many interesting and fundamental ideas and original
observations. In particular, they  studied the arithmetic of the
Fourier expansions $\sum a(n)w^n$ of noncongruence modular forms and
 observed that for some forms with algebraic coefficients the
sequences $\{a(n)\}$ have unbounded denominators. On the other hand,
it is known that if $f$ is a congruence holomorphic modular form,
the sequence $\{a(n)\}$ must admit bounded denominators
\cite{shim1}. This observation gives a simple and clear criterion
for distinguishing noncongruence modular forms.  Atkin and
Swinnerton  also found $p$-adic analogue of the Hecke three-term
recursion by demonstrating the coefficients $\{a(n)\}$ of
 some noncongruence modular forms satisfy
\begin{equation}\label{eq: 3term asd}
a(np)- A(p) a(n) + B(p)a(n/p)\equiv 0, \forall n\ge 1
\end{equation}
 modulo a power of $p$ depending on $\ord_p(n)$ for some constants
 $A(p)$ and $B(p)$. In \cite{a-sd}  a proof of the weight 2 dimension 1
case was provided. The modularity of this case would
 have followed from the Taniyama-Shimura-Weil conjecture.
   For a few examples of dimension 2 or more, Atkin and Swinnerton-Dyer found that such
 three-term congruences still exist with the forms
 diagonalized $p$-adically for each separate $p$ and the $A(p)$ being over algebraic
 number fields in general different for each $p$. Several people (\cite[etc]{ditters76,stienstra85})   studied Atkin and
 Swinnerton-Dyer (ASD)
congruences. The
 general situation for the weight 2 cases has been confirmed by several
  papers  \cite{hon2,car71,ditters76,Katz81}.
\smallskip

 In \cite{sch85b}, using deep algebraic geometry and representation theory, Scholl  associated each
 $d$-dimensional space of noncongruence cuspforms  with  $2d$-dimensional $l$-adic
 Scholl representations and  showed the existence of $(2d+1)$-term $p$-adic recursions under  general conditions.
 When $d=1$ Scholl  fully established the 3-term ASD
 congruences \eqref{eq: 3term asd}.   However when $d>1$, the existence of a fixed basis for each
 of which satisfies 3-term ASD congruences for all $p$ (referred to
 as simultaneous diagonalization)
 is rare.

Around the same time interesting connections between modular forms
for congruence and noncongruence subgroups began to merge. In
\cite{sch88,sch93} Scholl gave examples of 1 dimensional spaces of cuspform  with
 the $A(p)$
 being eigenvalues of congruence cuspforms by
 showing the  $l$-adic Scholl representations associated to those cases are
 ``modular". In \cite{berger00}, Berger found a quadratic
 relation satisfied by 2 noncongruence cuspforms and 2 congruence
 cuspforms.
\smallskip

 In
\cite{a-sd,sch88} noncongruence subgroups with small indices are
given considerable attention. Recently, investigation has been
focused on the so called noncongruence character groups, which are
``almost" congruence. These noncongruence subgroups are constructed
by assigning characters to the generators of congruence
subgroups. The Lattice groups considered by Rankin \cite{rankin67} are
examples of them. Klein and Fricke used ``Wurzelmodul" to describe a
Hauptmodul of a genus zero noncongruence character group defined by
a cyclic character. Here, we will associate the Fermat curve
$x^n+y^n=1$ with the modular curve for the intersection of two
particular character groups.

In \cite{lly05}, Li et al. gave a $2$-dimensional space of
 cuspforms for a noncongruence character group which
satisfies simultaneous diagonalization. Moreover, it was shown that
$A(p)$ are eigenvalues of explicit newforms. Later, Fang et al.
found a few more examples of this kind \cite{L5}. As a matter of
fact, simultaneous diagonalization can rarely be achieved. More
recently, Atkin et al \cite{all05} and Long \cite{long061}  fully
exhibited two cases which resemble the case in \cite{lly05} but are
slightly more general. In each case, a basis depending on the
congruence class of $p$ modulo a fixed integer is found, for each of
the basis function the coefficients satisfy a 3-term ASD congruence
relation. Meanwhile, modularity of each case was verified
individually. These modularity results rely heavily on the
decompositions of the associated $l$-adic Scholl representations or
their restrictions to suitable Galois subgroups. In contrast, Scholl
exhibited other cases where the $l$-adic Scholl
representations are not reducible \cite{sch04}.
\medskip

 This expository paper is organized
in the following order. Some background topics are introduced in
Section \ref{sec:background}. Section \ref{sec:diffences} is devoted
to some differences between congruence and noncongruence
subgroups. Modular forms and some related topics will be addressed
briefly in Section \ref{sec:mf}. In Section \ref{sec: unbounded
denominators}, we will show the Fourier coefficients of some
noncongruence modular forms have unbounded denominators. In Section
\ref{sec:modularity} we will address $l$-adic Scholl representations
attached to noncongruence cuspforms and related modularity results
briefly. Section \ref{sec:asd} is about the three-term
Atkin-Swinnerton-Dyer congruence relations. Throughout this paper,
several examples will be discussed and  all proofs presented  are
elementary. The author has not made an attempt to give a
comprehensive review of the vast literature. This survey mostly
reflects the author's perspectives on these topics based on her
recent research activities.
\medskip

This  paper has grown out of the author's talk delivered at the
Banff conference. The author is indebted to A.O.L. Atkin. The
introduction of this article is written partly based on the valuable
historic recounting and analysis of the development of noncongruence
modular forms by A.O.L. Atkin (communicated through private
correspondence). His wisdom and original ideas have sparked much of
the discussion here. It is the author's great pleasure  to thank
Wenching Winnie Li for her constant
 discussions and many insightful comments. Some results below are  cited
from some unpublished notes of Li and the author. The author would
also like to thank the conference organizers and in particular
Noriko Yui for inviting her to attend this wonderful conference and
giving her an opportunity to present this talk. She is also very
grateful to many people such as Ron Livn\'e, Siu-Hung Ng, Jan
Stienstra, Yifan Yang,  Don Zagier, for their valuable discussions
and comments. The author specially thanks Chris Kurth for his
remarks on an earlier version of this paper, Tonghai Yang for
pointing out several useful references, and the referee for many suggestions on improving the clarity of this paper.

\section{Backgrounds}\label{sec:background}
In this article, we use $p$ and $l$ for  prime numbers and $\Q_l$
for the field of $l$-adic numbers. We adapt mostly common notations.
We use $w$ for all local parameter at infinity instead of $q^{\mu}$
where $\mu$ varies in different cases.

\subsection{Finite index subgroups of $P\SL$}
It is well-known that $\SL$ is generated by two elements
$E=\begin{pmatrix}   0&1\\-1&0
\end{pmatrix}$ and $V=\begin{pmatrix}
  1&1\\-1&0
\end{pmatrix}$ \cite{Rankin-book-mf}. Identifying elements in $\SL$
which  differ by a sign, one obtains the modular group $P\SL$. It
is well-known that $P\SL$ is the free-product of two elements with
orders 2 and 3 respectively, namely the images of $E$ and $V$ in
$PSL_2(\Z)$
\cite{Rankin-book-mf}. In this article, we use $\Gamma$ to denote
finite index subgroups   of $\SL$. We also use $\G$ to denote
the image of $\pm\G/\pm I$ in $P\SL$ if there is no confusion.

Congruence subgroups are a class of subgroups which are
easy to describe. For example, given a positive integer $n$, the
following are some well-known congruence subgroups.
\begin{eqnarray*} \G^0(n)&=& \left \{ \gamma \in \SL,
\gamma=\begin{pmatrix}   *&0 \\ *&*
\end{pmatrix} \mod n \right \}.\\
\G^1(n)&=& \left\{ \gamma \in \SL, \gamma=\begin{pmatrix}   1&0 \\ *&1
\end{pmatrix} \mod n \right \}.\\
\G(n)&=&\left \{ \gamma \in \SL, \gamma=\begin{pmatrix}   1&0 \\ 0&1
\end{pmatrix} \mod n \right \}.
\end{eqnarray*}

 In \cite{Millington691,Millington692}, Millington
showed that up to isomorphism there is a one-to-one correspondence between finite index
subgroups of $P\SL$ and a pair of  permutations
 of orders 2 and 3 respectively which generate a finite transitive permutation group.

\begin{example}
 The group $\G^1(3)$ is an index 4 subgroup of $\SL$. It has 4 right coset
 representatives $\gamma_1=\begin{pmatrix}
   1&1\\0&1
 \end{pmatrix}, \gamma_2=\begin{pmatrix}
   1&2\\0&1
 \end{pmatrix}, \gamma_3=\begin{pmatrix}
   1&3\\0&1
 \end{pmatrix}, \gamma_4=\begin{pmatrix}
   0&-1\\1&0
 \end{pmatrix}$. For simplicity, we will use 1,2,3,4 to denote the corresponding cosets.   Multiplying $E$ on the right will permute the right
 cosets. More explicitly, $E=(12)(34)$ is an order 2 permutation.
 Similarly, $V$ corresponds to an order 3 permutation $V=(134)(2)$.
 Millington's theorem says $\G^1(3)$ corresponds uniquely (up to
 isomorphism) to the pair of permutations $(12)(34)$ and $(134)(2)$.
 Note that in this case $\begin{pmatrix}
   1&-1\\0&1
 \end{pmatrix}=VE=(123)(4)$. Since  the group of stabilizers of the cusp $\infty$ is generated by $\begin{pmatrix}
   1&-1\\0&1
 \end{pmatrix}$, the  corresponding
 cycle $(123)$  is called  the marked cycle in the
 disjoint cycle representation of $VE$ with respect to $\infty$.
\end{example}
\subsection{KFarey}
  The permutation representation of these groups is of great
  importance to the computational aspects of finite index groups of the modular group.  For example, in
  \cite{hsu96}, Hsu gave an algorithm for identifying congruence subgroups based on the permutation
  representations. Recently, Kurth has developed a SAGE package  called ``KFarey" which has implemented many useful functions for
  finite index
  subgroups such as Hsu's algorithm mentioned above, Farey symbols by
  Kulkarni \cite{Kulkani85}, determining whether a given element is
  in a group specified by a given permutation representation and a
  ``marked" cycle in the permutation representation of $VE$ corresponding to infinity,
  etc. The ``KFarey" can be downloaded from the following website
  \begin{quote}
    \emph{http://www.public.iastate.edu/~kurthc/research/index.html}
  \end{quote}

\subsection{Modular curves}\label{subsec:modcurve}
 The elements $\begin{pmatrix}
  a&b\\c&d
\end{pmatrix}\in \G$ act on the upper half plane $\fH= \left \{ z\in \C\,|\, \text{Im} z>0 \right \}$
via the linear fractional transformation:
$$\begin{pmatrix}
  a&b\\c&d
\end{pmatrix} \cdot z=\frac{az+b}{cz+d}, \quad  \forall z\in \fH
.$$  An element $\gamma=\begin{pmatrix}   a&b\\c&d
\end{pmatrix}\in \SL$ is said to be elliptic (resp. parabolic or hyperbolic) if the value of $|a+d|$ is $<2$
(resp. $=2$ or $>2$). Each elliptic (resp. parabolic) element has a
unique fixed point on $\fH$ (resp. $\Q \cup \{ \infty\}$).
Hyperbolic elements have no fixed points on
$\overline{\fH}=\fH\cup\Q\cup \{ i \infty\}$. (Whether a group has
elliptic elements or not can be easily read off from its permutation
representation. For example, $\G^1(3)$ has order 3 elliptic points
as the left coset labeled by 2 is fixed by $V$.)   The orbit space
$\fH/\G$ can be compactified canonically by adding a few isolated
points called the cusps of $\G$ (equivalence classes in
$(\{\infty\}\cup \Q)/\G$). The compactified curve, denoted by
$X_{\G}$ is called the \emph{modular curve } associated with $\G$. By
convention, the genus of $X_{\G}$ is also called the genus of $\G$.
For convenience, we will assume that all modular curves here are
defined over $\Q$ unless specified otherwise. When the curve
$X_{\G}$ has genus 0, the field of meromorphic functions on $X_{\G}$
is generated by a single element, which is called a Hauptmodul of
$\G$. For example, the classical modular $j$-function is a
Hauptmodul of $\SL$.

The importance of finite index subgroups of $\SL$ is seen in the
following theorem.
\begin{theorem}[Belyi \cite{belyi79}]
Any smooth compact complex projective curve defined over
$\overline{\Q}$ is isomorphic to the modular curve  of some finite
index subgroup of $\SL$.
\end{theorem}

On the other hand, one should know that the above correspondence is
not unique. Please refer to Section \ref{sec:Fermat} for a
connection between Fermat curves and modular curves which
demonstrates Belyi's theorem above.

Up to conjugation, the stabilizer of each cusp $c$ of $\G$ is
generated by a matrix of the form $\pm \begin{pmatrix} 1&m\\0&1
\end{pmatrix},$ where $m\ge 1$ is called the cusp width of $c$. The \emph{level}
of the group $\G$, denoted by $\text{lev}(\G)$, is  the least common
multiple of all cusp widths of $\G$. It was defined by Wohlfahrt
\cite{wohlfahrt64} which generalizes Klein's level definition for
congruence subgroups.

\subsection{Elliptic modular surface}\label{sec:ems} An elliptic
surface $$\pi: S \rightarrow C$$ is a two dimensional complex
variety $S$ together with a morphism map $\pi$ to the base curve $C$
such that for almost all $t\in C$, $\pi^{-1}(t)$ is an elliptic
curve. When $\pi^{-1}(c)$ is not an elliptic curve, it is called a
special fiber. Special fibers of the N\'eron models of elliptic
surfaces have been classified by Kodaira \cite{kod1} and N\'eron
\cite{ner1}. An elliptic surface $\pi: S \rightarrow C$ is called
semistable if all of its special fibers are of $I_n$ type (in
Kodaira notation) and the list of the (un-ordered) numbers
$[n_1,n_2,\cdots,n_r]$ in the subscripts of the special fibers is
called the configuration of the semistable elliptic surface.

Let $\C$ be the ground field. Given an elliptic fibration $\pi:
S\rightarrow C$, let $\Sigma$ be the finite collection of points on
$C$ corresponding to special fibers. Via a standard procedure (cf.
\cite{kod1}) one obtains the monodromy representation of the
fundamental group $\pi_1(C\setminus \Sigma, t_0):$
$$\rho_{\pi} : \pi_1(C\setminus \Sigma, t_0) \rightarrow \GL_2(\C),$$ where
$t_0$ is a fixed reference point on $C\setminus \Sigma$. The image
of  $\rho_{\pi}$, well-defined up to a conjugation corresponding to
the difference choices of $t_0$, is called the monodromy group of
$\pi$. Nori gave criterions for when the monodromy group of an
elliptic surface $\pi: S \rightarrow C$ gives rise to a finite index subgroup
of the modular group \cite{nor1}. By Nori's theorem, the monodromy group of a
semistable elliptic surface with maximal possible Picard number is a
finite index subgroup of $\SL$. Here for a surface $S$, the
$\Z$-rank of its group of divisors modulo linear equivalence is
called the Picard number of $S$.

 In \cite{shio1}, for every finite index subgroup $\G$ of $\SL$ not containing $-I_2$, Shioda
constructed an elliptic modular surface $\pi: \E_{\G} \rightarrow
X_{\G}$ associated with $\G$. In particular, for an index 24 genus
zero subgroup $\G$,  $\E_{\G}$ is an elliptic K3 surface. K3
surfaces are simply connected compact complex surfaces with trivial
canonical bundles. They are two dimensional Calabi-Yau manifolds and
occupy an important role in the classification of compact complex
surfaces (cf. \cite{bpv}).   Miranda and Persson have classified
all configurations of semistable elliptic K3 surfaces with Picard
number 20 \cite{m-p2}.
 Each K3 surface of this kind
corresponds to an index 24 subgroup with 6 cusps and its cusp widths
are given by the configuration of the surface.  (In \cite{m-b04},
Montanus and Beukers  computed algebraic models of those semistable
elliptic K3 surfaces defined over $\Q$.) Among the 112 possible
extremal semistable K3 configurations, 9 cases correspond to
congruence subgroups. In \cite{ty2}, Topp and Yui have mentioned
that among the remaining 103 configurations corresponding to
noncongruence subgroups, 10 can be obtained via double covers of the
semistable extremal rational elliptic surfaces classified
 by Beauville \cite{bea1} ramifying at only two points. For example,
  a semistable K3
with a configuration  $[1,1,1,1,10,10]$ is a double cover of a
semistable  rational surface with a configuration $[1,1,5,5]$ (cf.
Example \ref{eg:G_n} below).

\subsection{Most finite index subgroups are noncongruence}

\noindent Rademacher conjectured that there are only finitely many genus zero
congruence subgroups (c.f \cite{Knopp-Newman65}). In
\cite{Dennin751}, Dennin showed  that there are only finitely many
congruence subgroups of a given genus. In \cite{seb1}, Sebber
classified all
 elliptic elements free   genus zero congruence subgroups. There are 33 such groups in total.
 In \cite{CLY04}, Chua, Lang, and Yang gave explicit descriptions for all the genus zero congruence subgroups. In \cite{cummins-pauli03}, Cummins and Pauli classified
 all congruence subgroups with genus less than 25.

 On the other hand,   Newman
 showed that there are infinitely many genus zero
 noncongruence subgroups  \cite{newman65}. In
\cite{a-sd}, Atkin and Swinnerton-Dyer  pointed out that
noncongruence subgroups predominate congruence subgroups.  In
\cite{jon79}, Jones proved that
 there are infinitely many noncongruence
subgroups of a given genus (cf. also  \cite{jon86}). Let $N_c(n)$
(resp. $N_{nc}(n)$) be the number  of congruence (resp.
noncongruence) subgroups of index $n$ in $P\SL$. Stothers gave an
upper bound for $N_c(n)$ and showed that
$N_c(n)/N_{nc}(n)\rightarrow 0$ when $n \rightarrow \infty$
\cite{Stothers84}.

\subsection{Noncongruence character groups}\label{subsec:ncc}
It is well-known that for any positive integer $n>2$, the finite
index subgroup $\G^1(n)$ of $\SL$ is a normal subgroup of $\G^0(n)$
with the quotient group $\G^0(n)/\G^1(n)$ isomorphic to
$(\Z/n\Z)^{\times}$. Such a relation provides a canonical
decomposition for the spaces of modular forms (cf. Section
\ref{sec:mf}) for $\G^1(n)$ in terms of the characters of
$(\Z/n\Z)^{\times}$ (cf. \cite{kob1}). In the same spirit,
character groups are defined as follows:
\begin{defn}\label{defn:ncc}
A group $\G$ is called a \emph{character group}  of a finite index
subgroup $\G^0$ if $\Gamma$ is normal in $\G^0$ with abelian
 quotient. Equivalently, a
character group is the kernel of an onto homomorphism
\begin{equation}
  \varphi: \G^0 \rightarrow G
\end{equation} where $G$ is a finite abelian group written multiplicatively.
\end{defn}

By the Fundamental Theorem of Finite Abelian Groups, it suffices to
consider when $\G^0/\G$ is cyclic.
\medskip

When working with character groups, we  distinguish the
following two types of character groups. 

\begin{defn}
  A character group $\G=\ker \varphi$ is said to be of type I if $\varphi(\gamma)\neq
  1$ for some parabolic $ \gamma \in \G^0$; otherwise, it is said to be of
  type II.  $\G$ is said to be of type II(A), if the images of all
  elliptic and parabolic elements of $\G^0$ are 1 under $\varphi$.
\end{defn}

Note that if $\G$  is a type II character group of $\G^0$ then the levels of $\G$ and $\G^0$ are the same; this is not true for
type I character groups in general. It has been shown that almost
all type II character groups are noncongruence \cite{Kurth-Long06}.
If $\G$ is a type II(A) character group of $\G^0$, then the modular
curve for $\G$ is a strictly $[\G^0:\G]$ to 1  unramified covering
map of the modular curve for $\G^0$. By definition, the image of any
elliptic element $\gamma$ in $\G^0$ under $\varphi$ has to divide
the order of $\gamma$ which is a finite positive integer.

\medskip

In \cite{lly05, all05, long061}, the following family of type I
 character groups is considered.
\begin{example}\label{eg:G_n}
Let $$\G^1(5)=\left \{ \gamma\in \SL \big| \gamma \equiv \begin{pmatrix} 1&0\\
* &1
\end{pmatrix} \mod 5 \right \}.$$ This is an index 12 elliptic element free congruence subgroup of
$\SL$ with 4 cusps: $ 0,-5/2, -2,\infty$. The surface $\E_{\G^1(5)}$
is semistable with configuration  $[1,1,5,5]$.  It is one of the 6
semistable extremal rational surfaces classified by Beauville
\cite{bea1}. 
This group is generated by 4 elements denoted by $\gamma_{\infty},
\gamma_{-2}, \gamma_0,\gamma_{-5/2}$ (each stabilizes one cusp as
indicated by the subscripts)  subject to one relation
$$\gamma_{\infty}\gamma_{-2}\gamma_{0}\gamma_{-5/2}=I_2.$$ Let
$\varphi_n$ be the homomorphism from $\G^1(5)$ to $\C^{\times}$
sending $\gamma_{\infty}$ to $\omega_n$, $\gamma_{-2}$ to
$\omega_n^{-1}$, and $\gamma_0, \gamma_{-5/2}$ both to 1, where
$\omega_n=e^{2\pi i/n}$. The kernel of $\varphi_n$,  denoted  by
$\G_n$, is a  type I character group of $\G^1(5)$. We list a few
properties of $\G_n$ here:
\begin{itemize}

\item[1.]    $\G_n$ is  noncongruence
 when $n\neq 1,5$ (cf. Section \ref{sec: unbounded
denominators}). \item[2.] $\G_n$ is elliptic element free with index
$12n$ in $\SL$ and has $2+2n$ cusps.
\item[3.]  $\G_5$ is the last group in Sebbar's list \cite{seb1}.
\item[4.] The elliptic modular surface for $\G_2$ is a K3 surface with configuration $[1,1,1,1,10,10]$.
\end{itemize}
\begin{figure}
\begin{center}
\scalebox
{0.3} 
{
\includegraphics*{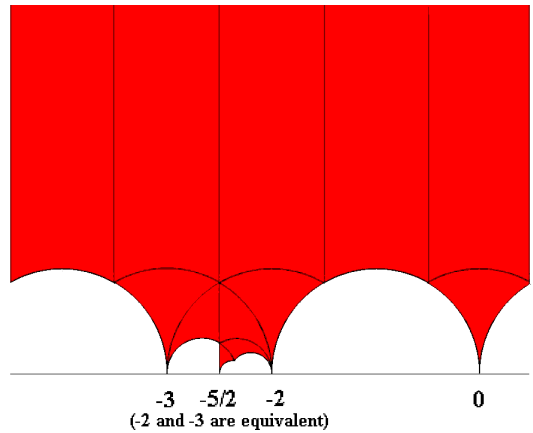}
}
\end{center}
\caption{Fundamental domain for $\G^1(5)$, drawn by the Fundamental
Domain Drawer of
 Verrill}
\end{figure}
\end{example}

\begin{example}
The group $\G^0(11)$ is genus 1, torsion free, and index 12 in
$\SL$. It has two cusps: $\infty$ and $0$. It is generated by 4
generators: $g_{\infty}$ and $g_0$ which stabilize $\infty$  and 0
respectively and two parabolic generators $A$ and $B$ subject to one
condition
$$g_{\infty}g_0ABA^{-1}B^{-1}=I_2.$$ Then the kernel of any homomorphism which sends
$g_{\infty}$ and $g_0$ to the identity and $A$ and $B$ to roots of
unity is a type II character of $\G^0(11)$.  Among all such groups,
there is only one, namely $\G^1(11)$, which is congruence. By the
Hurwitz formula, the genus of such a group remains  1 and as we have
remarked before the level of this group is 11.  For any type II
character group $\G$ of
  $\G^0(11)$, its modular curve $X_{\G}$ is an elliptic curve which
  is
  isogenous to $X_{\G^0(11)}$.
\end{example}

\section{Some distinctions between congruence and noncongruence
subgroups}\label{sec:diffences}\subsection{}One fundamental
difference between congruence and noncongruence subgroups is that
elements of the congruence subgroups can be described by congruence
relations while the elements in the noncongruence subgroups cannot.

\subsection{}Conjecturally, Fourier coefficients of modular forms for these two
types of groups also behave differently. Please refer to Section
\ref{sec: unbounded denominators} for a more detailed discussion.

\subsection{} The following is a special case of the finiteness theorem of Thompson \cite{Thompson80}.
Let $\G$ be a finite index subgroup of the modular group
 and $\mathcal K$ the set of congruence subgroups of
$\G$, then
\begin{equation}
  \lim_{K\in \mathcal K} \frac{g(K)}{|\G: K|}=\frac{1}{2} \chi(\G),
\end{equation}in the sense that for any $\varepsilon>0$, there are only finitely many elements $K\in \mathcal K$ such that
$$\left | \frac{g(K)}{|\G:K|}-\frac{1}{2}\chi(\G)\right |>\varepsilon.$$
Here $g(\G)$ is the genus of $\G$,
\begin{equation}
  \chi(\G)=2(g(\G)-1)+t+\sum_{i=1}^m\left (1-\frac{1}{e_i}\right ),
\end{equation} is the Euler characteristic of $\G$, $t$ is the
number of cusps of $\G$, $m$ is the number of elliptic points of
$\G$, and the $e_i$'s are the orders of the elliptic points in
$P\SL$.

In comparison,
 there are infinitely many noncongruence
subgroups of a given genus \cite{jon79}. Thus,
 the finiteness theorem of Thompson \cite{Thompson80} does not hold if the congruence condition on the set of elements in $\mathcal K$ is removed.

\subsection{}
Let $A(\fH/\G)$ denote the hyperbolic area of a fundamental domain
of $\G$. Let $$\Delta=-4(\text{Im} z)^2\frac{\partial^2}{\partial z
\partial \bar{z}}$$ be the Laplace operator of the Poincar\'e metric on
$\fH$. It admits a self-adjoint non-negative extension in the space
of square integrable functions on $\fH/\G$, which is denoted by
$\Delta_{\G}$. Let $\l_1$ be the smallest positive  discrete
eigenvalue of $\Delta_{\G}$ (if it exists).
\begin{theorem}[Selberg \cite{sel65}]
If $\G$ is a congruence  subgroup, then  $\l_1\ge 3/16$.
\end{theorem}
\begin{theorem}[Zograf \cite{Zograf91}]\label{thm:Zograf}Let $\G$ be a discrete subgroup of
$PSL_2(\mathbb R)$ with hyperbolic area $A(\fH/\G)<\infty$, and
assume $A(\fH/\G)\ge 32 \pi\cdot (g(\G)+1)$. Then $\l_1$ exists and
\begin{equation}\label{eq:Z}
  \l_1< \frac{8\pi\cdot (g(\G)+1)}{A(\fH/\G)}.
\end{equation}
\end{theorem}

\begin{theorem}
  For any $\varepsilon>0$, there exists a noncongruence
  subgroup with $0<\l_1<\varepsilon$.
\end{theorem}
\begin{proof}
  As there are genus 0 noncongruence subgroups with $A(\fH/\G)$
  arbitrarily large (cf. Example \ref{eg:G_n}), the claim follows from \eqref{eq:Z}.
\end{proof}

\section{Modular forms}\label{sec:mf}\subsection{} A holomorphic (or meromorphic) \emph{modular form}
$f$ for  $\G$ is a holomorphic (or
meromorphic) function defined on $\overline{\fH}$ satisfying
 \begin{equation}\label{eq:mf}
f(\frac{az+b}{cz+d})=(cz+d)^{k}f(z), \quad \forall z \in \fH,\,
\begin{pmatrix}   a&b\\c&d
\end{pmatrix}\in \G.
\end{equation} The
number $k$ is called the weight of the modular form $f$.  A weight
$0$ meromorphic modular form is called a \emph{modular function }for $\G$.
All modular functions for $\G$ together form a field, denoted by
$\frak M_{\G}$, which is transcendental over $\C$. When $\G$ has
genus zero, $\frak M_{\G}$ is generated by any Hauptmodul for $\G$.

A holomorphic modular form is called a cuspform if it vanishes at
every cusp of $\G$. Let $S_k(\G)$ denote the space of holomorphic
weight $k$ cuspforms for $\G$.  It is a finite dimensional vector
space over $\C$. Explicit formulae for the dimension of $ S_k(\G)$
are available in \cite{shim1}. Cuspforms are important invariants of
$\G$ and $X_{\G}$. For instance weight 2 cuspforms for $\G$ are
invariant differential 1-forms of $X_{\G}$. A cuspform of integral
weight larger than 2 is a differential form on $X_{\G}$ with
coefficients in a line bundle.

For any $\gamma=\begin{pmatrix}   a&b\\c&d
\end{pmatrix}\in \SL$, define the stroke operator $|_{\gamma}$ action on a
weight $k$ modular form $f$ to be
\begin{equation*}
  f|_{\gamma}=(cz+d)^{-k}f(\frac{az+b}{cz+d}).
\end{equation*}

Theta functions and Eisenstein series constitute important classes
of modular forms. They are described in most books on modular forms
(\cite[etc]{ogg,shim1,kob1,d-s-mfbook}). In this paper, we use $w$
to denote $e^{2\pi i z/\mu }$ where $\mu$ is the width of the cusp
at infinity. Accordingly, any modular form has an expansion of the
form $\sum_{n\ge n_0} a(n) w^n$.

\begin{example}\label{eg:E12}
There are two particular holomorphic weight 3 Eisenstein series for
$\G^1(5)$: $E_1$ and $E_2$, which vanish at all the cusps of
$\G^1(5)$ except $\infty$ and $-2$ respectively. Both $E_1$ and
$E_2$ have integral Fourier expansions at $\infty$. Moreover,
$E_1\cdot E_2=\eta(5z)^3\cdot\eta(z)^9$, here $\eta(z)=q^{1/24}\prod_{n\ge 1}(1-q^n), q=e^{2\pi i z}$ is the
Dedekind eta function. One interesting expression for $E_1$ is
$$E_1=1-\sum _{n\ge 1}b_n \frac{n^2w^n}{1-w^n}, $$ where
$b_n=2,1,-1,-2,0,2,1,-1,-2,0\cdots$  \cite{stienstra05}.

A
Hauptmodul of $\G^1(5)$ is
\begin{equation}\label{eq:t}
  t=\frac{E_1}{E_2}=w^{-1}+5+10w+5w^2-15w^3-24w^4+15w^5+\cdots,
  \quad
  w=e^{2\pi i z/5}.
\end{equation} Another  formula for $t$ is (cf. \cite{sebbar02})
$$t=\frac{1}{w} \prod_{n\ge 1} (1-w^n)^{-5 (\frac{n}{5})},$$ where
$\left (\frac{\cdot }{5}\right )$ is the Legendre symbol.

A Hauptmodul for $\G_n$ is
\begin{equation}\label{eq:tn} t_n=\sqrt[n]{t}.
\end{equation}
\end{example} In particular, it is  also known that $t_5$ is
a Hauptmodul for $\G(5)$ (when $w=e^{2\pi iz}$).

\subsection{Cuspforms for noncongruence character groups}
Due to the definition of noncongruence character groups, modular
forms for such a group $\G$ are closely related to modular forms for
$\G^0$. The following example illustrates such a relation.  When
$k>3$ we have the following result.
\begin{lemma}\label{eg:mfGn}Let $\G_n$ be as  in Example \ref{eg:G_n}.  Let $k>3$ and  $f_{n,k}\in S_k(\G^1(5))$ with simple
zeros at the cusps $-2,0,-5/2$ and an order $k-3$ zero at infinity
(and no other zeros). Then $$f_{n,k}\cdot t_n^j, j=1-n,\cdots,
(2k-3)n-1.$$ form a basis of $S_{k}(\G_n)$.

\end{lemma}
\begin{proof}  Using the dimension formula  in \cite{shim1}, we have $\dim
S_{k}(\G_n)= 2n(k-1)-1$. By the Riemann-Roch theorem, such a
cuspform $f_{n,k}$ exists and for any $j\in \{1-n,\cdots,
(2k-3)n-1\},$ $f_{n,k}\cdot t_n^j \in S_{2k}(\G_n).$ As these
functions are linearly independent, they form a basis.
\end{proof}
\subsection{Hecke operators}\label{sec:major difficulty}  Many good properties of congruence modular forms are
derived from Hecke operators. The major difficulty in understanding
modular forms for noncongruence subgroups is  the lack of
satisfactory Hecke operators. For a finite index subgroup $\G$,
define the $p$th Hecke operator $T_{p,\G}$ associated to $\G$ as
follows: Let $\a_p=\begin{pmatrix} p&0\\0&1
\end{pmatrix}.$ We decompose the set $\G \a_p\G=\{\gamma \a_p \gamma'| \g, \g'\in \G\}$ into a finite disjoint union $\bigcup_i \G  \a_i$. For any
modular form $f$  for $\G$, define \begin{equation*}   T_{p,\G} f=
\sum f|_{\a_i}.
\end{equation*} It is conjectured
by Atkin and has been shown by Serre (for $\G$ normal in $\SL$
\cite{tho89}), Scholl (for certain cases \cite{sch97}), and Berger
(in general  \cite{berger94}) that
\begin{theorem} For any prime $p\nmid \text{lev}(\G)$ and all $f\in S_k(\G)$
$$T_{p,\G}f=T_{p,\G^0} \circ \Tr_{\G}^{\G^0}f,$$ where
$\Tr_{\G}^{\G^0}: S_k(\G) \rightarrow S_k(\G^0)$ is the trace map.
\end{theorem}
Consequently for all $f\in  \, S_k(\G) \setminus S_k(\G^0)$,
$T_{p,\G}f=0$.

\section{Unbounded denominators}\label{sec: unbounded
denominators} \subsection{} \label{subsec:ubd}In what follows, we only consider
modular forms $\sum_{n\ge n_0} a(n)w^n$ with coefficients $a(n)$ in
a fixed number field.  Due to  Hecke operators, each space of holomorphic
congruence modular forms has a basis consisting of forms with
integral coefficients \cite{shim1}. Consequently, for every
congruence holomorphic modular form  with algebraic coefficients, the sequence
$\{a(n)\}$ has bounded denominators in the sense that there exists an algebraic number $M$ such that $M\cdot a(n)$ is algebraic integral for all $n$. In \cite{a-sd}, Atkin and
Swinnerton  observed that this is no longer the case for
noncongruence modular forms. Therefore, the sequence $\{a(n)\}$
having unbounded denominators  clearly implies $f=\sum _{n\ge n_0}
a(n)w^n$ is noncongruence.  At present, it remains unknown whether
$\{a(n)\}$ having bounded denominators will imply $\sum a(n)w^n$ to
be congruence. We will formulate a more precise conjecture as
follows:

Let ({UBD}) refer to the following property:
\begin{quote}
  \emph{Let $f$ be an arbitrary holomorphic integral weight modular form  for  $\G$  but not for $\G^0$ with algebraic Fourier
  coefficients at infinity. Then the Fourier coefficients of $f$ have
  unbounded denominators.}
\end{quote}
An open question is:
\begin{quote}
\emph{Does every noncongruence  subgroup satisfy the
condition (UBD)}?
\end{quote}

\begin{remark}
  In \cite{a-sd} and \cite{sch85b},   the coefficients $a(n)$  are normalized
  so that they are contained in
smaller fields. In the following discussion, we have avoided the
normalization on purpose for the following two reasons. First, such
a procedure is unnecessary for all cases below. Second, it is more
appropriate to address the bounded denominator property  prior to
the normalization as it is pointed out by D. Zagier that all forms
with algebraic coefficients can be normalized to have  integral
coefficients.
\end{remark}
\subsection{}

For any $n\ge 1$ and with a principal branch fixed, we formally
write
\begin{equation}\label{eq:nthroot}
 (1+x)^{1/n}= \sqrt[n]{1+x}=\sum_{m\ge 0}\frac{\left (\frac{1}{n} \right )_m}{m!}x^m,
\end{equation} where
$(\frac{1}{n})_m=\frac{1}{n}(\frac{1}{n}-1)\cdots
(\frac{1}{n}-m+1)$. The next lemma provides a simple but useful
criterion for detecting the unbounded denominator property in simple
cases.
\begin{lemma}\label{lem:unbounded}
 Let $n$ be any natural
 number and $f=1+\sum_{m\ge 1}a(m)w^m, a(m)\in \Z$ for all $m$.  In terms of (\ref{eq:nthroot}), we expand
$\sqrt[n]{f}=\sum_{m\ge 0} b(m)w^m, b(m)\in \Z[1/n]$. Let $p$ be a
prime factor of $n$. If there exists one
 $b(m)$ which is not $p$-integral, then
$$\limsup_{m\rightarrow \infty}-\ord_p(b(m)) \rightarrow \infty.$$
In other words, $\{b(m)\}$ has unbounded denominators.
\end{lemma}
\begin{proof}
  It is a special case of \cite[Lemma 2.7]{Kurth-Long06}.
\end{proof}
\subsection{}

\begin{example}\label{eg:quintic}
In \cite{COR-VI},  Candelas, de la Ossa, and Rodriguez-Villegas have
studied the  family of quintic Calabi-Yau 3-fold in $\mathbb P^4$
given by
$$\sum _{i=1}^5 x_i^5-5 \psi x_1x_2x_3x_4x_5=0.$$ They have shown that
$$\frac{1}{(5\psi)^5}=w+154w^2+179139w^3+\cdots\in \Z[w]], w=e^{2\pi  i t},$$
where $$t=\frac{1}{2\pi i} \frac{ \bar{\omega}_1
(\psi)}{\bar{\omega}_0 (\psi)}$$ is the a ratio of two particular
solutions of the Picard-Fuchs equations arising from the quintic
family. By the above lemma, the Fourier coefficients of $\psi$ have
unbounded denominators.
\end{example}
\subsection{}
\begin{example}Let $\G_n$ be as in  Example \ref{eg:G_n}. A Hauptmodul of
$\G_n$ is $t_n=\sqrt[n]{t}$ where the expansion of $t$ is as in
\eqref{eq:t}. By Lemma \ref{lem:unbounded}, when $n\neq 1,5$, the
coefficients of $t_n$ have unbounded denominators. Therefore $\G_n$
is noncongruence.
\end{example}

\begin{lemma}
  Let $f_{n,k}$ be the weight $k$ cuspform for $\G^1(5)$ as in Lemma \ref{eg:mfGn}, then the Fourier coefficients of
  $1/f_{n,k}$ have bounded denominators.
\end{lemma}
\begin{proof}
By our choice of $f_{n,k}$ whose poles are supported at the cusps, we know not only
   $f_{n,k}$  satisfies the bounded denominator property. Take any integral
   weight normalized newform $g$ for $\G^1(5)$. Clearly for any positive integer $m$,  both $g^m$ and $g^{-m}$ satisfy the bounded denominator property.
   When $m$ is large enough, $g^m(f_{n,k})^{-1}$ is also a holomorphic cuspform $h$ for $\G^1(5)$, hence $(f_{n,k})^{-1}=hg^{-m}$ also
satisfies the bounded denominator property. \end{proof}

\begin{prop}\label{prop:1}
Let $f$ be a weight $k>3$ cuspform for some $\G_n$ with coefficients
in some algebraic number field. If it is not a form for $\G_1$ or
$\G_5$, then the coefficients of
  $f$ at $\infty$  have unbounded denominators.
\end{prop}
\begin{proof}
  By Lemma \ref{eg:mfGn}, $f=f_{n,k}\cdot (\sum _{j=1}^na_jt_n^j)$ for some
  $a_j \in \C$.  Let $F[x]\in \C[x]$ be the nonconstant polynomial with minimal degree  such
  that the coefficients of $F(t_n)$ have bounded denominators. This is
  $F[x]=x^n$ if $5\nmid n$ or $F[x]=x^{n/5}$ otherwise. In the first
  case, if $f\notin S_k(\G^1(5))$, then $a_j\neq 0$ for some $j$ not
  divisible by $n$, hence $(\sum a_jx^j)$ is not a polynomial in
  $x^n$, i.e. the Fourier coefficients of $(\sum a_jt_n^j)$ have
  unbounded denominators. If $f$ satisfies the bounded denominator property so does $f(f_{n,k})^{-1}=(\sum a_jt_n^j)$
   which contradicts our assumption above.
      The second case can be handled in a similar manner.
\end{proof}

\subsection{The Fermat curves and modular curves}\label{sec:Fermat} In
the next example, we give one of the realizations of the Fermat
curves in terms of modular curves. It is motivated by a question
raised by Koblitz (through private conversation) regarding the
 relations between the Fermat curves and noncongruence subgroups.

  For any positive integer $n\ge 2$, $$x^n+y^n=1$$ is the Fermat curve
  and its properties have been studied extensively (cf.
  \cite[etc]{Gross-Rohrlich78,Kobliz-Rohrlich78,yui80,YangTH96}). Now let $G_1$ and  $G_2$ be two genus
  0  subgroups with Hauptmoduls $$x=\sqrt[n]{\l}, \quad y=\sqrt[n]{1-\l}$$ respectively, where the lambda function
  \begin{eqnarray*}
  \l(q)&=&16q\prod_{n\ge 1} [(1+q^{2n})^2(1-q^{2n-1})]^8 \\
  &=&16q(1-8q+44q^2-192q^3+718q^4-2400q^5+7352q^6+\cdots)\\
\end{eqnarray*} is a Hauptmodul for $\G(2)$.   (Not just  any arbitrary algebraic function of $\l$
  is another modular function. These two functions can be
  Hauptmoduls due to the special locations of the zeros and poles of $\l$ and $1-\l$).
  Clearly $x$ and $y$ satisfy $x^n+y^n=1$. It is not difficult to see that the Fermat curve above is bi-holomorphic to the modular curve for
  $\Phi(n)=G_1\cap G_2$. In the literature, $\Phi(n)$ is the so-called Fermat group. By Lemma
  \ref{lem:unbounded}, it is easy to check that when $n\neq 1,2,4,8$ the Fourier coefficients of the expansions of $x$ and $y$ in $q$ at infinity will
  have unbounded denominators, which implies  $\Phi(n)$ is noncongruence. By using the idea used in the proof of
  Proposition \ref{prop:1}, we can show that when $n\neq 1,2,4,8$, the Fermat group $\Phi(n)$ satisfies the
  condition (UBD) (cf. Section \ref{subsec:ubd}).

\subsection{}
Recently,
 Kurth has verified that for each of
 the following K3 configurations \cite{Kurth-K3}
$$[1,1,1,1,2,18],[1,1,1,1,10,10],[1,1,2,2,6,12],[1,1,2,2,9,9],[1,1,2,5,5,10]$$
$$[1,1,4,6,6,6],[2,2,2,3,3,12],[2,2,5,5,5,5],[2,3,3,4,6,6],[3,3,3,3,6,6]$$

\noindent arising as double covers of rational surfaces mentioned in
subsection \ref{sec:ems} and  \cite{ty2}, the
  Fourier coefficients of every possible Hauptmodul of these
  monodromy groups
 have unbounded denominators.

\subsection{}In general, it is relatively easy to construct modular
forms for type I character groups whose corresponding homomorphisms
send all elliptic and hyperbolic generators to 1 (cf. Example
\ref{eg:G_n} and Lemma \ref{eg:mfGn}). To handle all type II or
II(A) character groups of a fixed congruence group is more
challenging.  The following theorem gives a partial positive answer to
the open question mentioned at the beginning of this section for
type II(A) character groups.
\begin{theorem}[\cite{Kurth-Long06}]
  Let $\G^0$ be any genus 1  congruence   subgroup whose modular curve is  defined over a number field. If there
  exists a prime  $p$ such that every index-$p$ type II(A) character group
  of $\G^0$  satisfies (UBD), then  there exists $ c(\G^0)>0$  such  that for any $X\gg 0$,
\begin{equation}
\#\left \{\text{subgroup } {\G} \text{ of  } \G^0 \ | \
[\G^0:\G]<X,  \text{ and } \G \text{ satisfies (UBD)}
\right \}>c(\G^0)\cdot X^2.
\end{equation}
\end{theorem}It is computationally feasible to check the conditions above.
For example $\G^0=\G^0(11)$ is considered in  \cite{Kurth-Long06}
via some results of Atkin in \cite{Atkin67}.

\section{Modularity}\label{sec:modularity}

\subsection{$l$-adic Scholl representations}Let  $\G$ be a noncongruence  subgroup with $X_{\G}$ defined over $\Q$ and
 $k>2$ be an integer. Let $d=\dim S_k(\G)$.  For every prime $l$, Scholl has constructed a $2d$-dimensional
$\Q_l$-vector space $W$ from weight $k$ cuspforms for $\G$ such that
there is a continuous homomorphism
$$\rho: \Gal(\overline{\Q}/\Q)\rightarrow \text{Aut}(W),$$ (cf. \cite{sch85b}). These
representations satisfy good properties: They  ramify only at
 a finite set of primes.
For each prime $p$ such that $\rho$ is unramified, the
characteristic polynomial of the Frobenius element at $p$ has
integral coefficients. It has been verified that Scholl's results
also hold for odd weight cases (cf. \cite{L5}). In \cite{sch88},
Scholl  investigated the 2-dimensional $l$-adic representation
attached to the unique weight 4 genuine cuspform for an index 9
noncongruence subgroup. He showed that it is isomorphic to the
$l$-adic representation attached to a weight 4 level 14 congruence
newform. Scholl has given a few more unpublished examples of this
kind \cite{sch93}. 

\subsection{} In \cite{b-s}, Stienstra and Beukers  investigated several K3
surfaces with Picard number 20 over $\C$. These K3 surfaces are
constructed as double covers of  $\C P^2$ with branching points
along degree 6 polynomials. According to a result of Shioda and
Inose \cite{s-i}, the L-function of each algebraic K3 surface with
Picard number 20 mainly consists of the L-function of a weight 3
congruence Hecke eigenform. In other words, the $l$-adic
representations attached the nontrivial part of the $l$-adic second
\'etale cohomology of  such a K3 surface are modular. Stienstra and
Beukers have identified the Hecke newforms corresponding to each K3
surface in their consideration and have constructed elliptic fibers
for those K3 surfaces including the one with configuration
$[1,1,1,1,10,10]$. However, using Sebbar's classification
\cite{seb1}, one sees several monodromy groups attached to these
elliptic fibrations are noncongruence.  Combining the above
discussions with the theory of elliptic modular surfaces by Shioda
\cite{shio1}, we have
\begin{theorem} Let $\G$ be a finite index subgroup of the modular group which is conjugate to a
monodromy group of
 a K3 surface over $\Q$ with Picard number 20. Then the $l$-adic Scholl representation attached to
$S_3(\G)$ is isomorphic to the $l$-adic representation attached to
certain Hecke eigenform.
\end{theorem}
This modularity result has significant influence on \cite{lly05} and
it explains why weight 3 cuspforms are given more attention. As a
special case of the above theorem, we know the $l$-adic Scholl
representation attached to $S_3(\G_2)$ is modular.
\smallskip

\subsection{}
Starting from \cite{lly05}, investigations have been focused on
 higher dimensional $l$-adic Scholl
representations. In particular, $l$-adic Scholl representations
attached to $S_3(\G_n)$ have been examined extensively. To study
these representations, it is natural to decompose the representation
spaces into subrepresentation spaces by using available operators.
Normalizers of $\G$ in $SL_2(\mathbb R)$ become important as they give rise to
involutions on the corresponding $l$-adic representation spaces. Use
$\G^{N}$ to denote the group of normalizers of $\G$ in $SL_2(\mathbb R)$.
Clearly $\G^0\subset\G^{ N}$. (M.L. Lang gave an algorithm for
finding $\G^{N}$ for a given $\G$ \cite{Lang02normofmg}). Moreover,
$\G^{N}/\G$ acts on the finite dimensional spaces $S_k(\G)$. We
assume the action of any element in $\G^N/\G$ is defined over a
number field $F$. Consequently, the $l$-adic School representation
above, when restricted to the subgroup $\Gal(\overline{\Q}/F)$,  can
be decomposed into a direct sum of subrepresentations. In some ideal
cases, these subrepresentations are 2-dimensional and are shown to
be isomorphic to some $l$-adic Deligne representations attached to
congruence newforms \cite[etc]{lly05,L5,all05,long061}.

\subsection{}With the new exciting progress on major conjectures
like the Fontaine-Mazur conjecture \cite{Kisin06-FMconjecture}, and
the Serre conjecture \cite{KW1, Kisin07-Serre}, it is promising that one can establish the
modularity for a large class of $l$-adic Scholl representations.

\section{Atkin and Swinnerton-Dyer congruences}\label{sec:asd}
\subsection{}When $\G$ is a congruence subgroup, we have the following
classical result due to Hecke and Petersson. For any positive
integer $k$, the space $S_k(\G)$ has a basis consisting of Hecke
eigenforms. Namely, for each of the basis element $\sum_{n\ge 1}
a(n) w^n$ there exists a suitable Dirichlet character $\chi$ such
that for any prime $p$
\begin{equation}\label{eq:hecke3term}
a(np)-a(p)a(n)+\chi(p)p^{k-1}a(n/p)=0, \quad \forall \ n\in \mathbb N.
\end{equation} However, when $\G$ is a noncongruence subgroup, the Hecke operator theory
is unsatisfactory as mentioned in Section \ref{sec:major
difficulty}. Consequently, such a result no longer holds for
noncongruence cuspforms. While working with explicit examples, Atkin
and Swinnerton-Dyer have found that three-term recursions like
\eqref{eq:hecke3term} still hold in a $p$-adic sense. What Atkin and
Swinnerton-Dyer have observed is the following: given $S_k(\G)$,
then for almost all primes $p$, the space $S_k(\G)$ has a basis (pending on $p$)
consisting of forms of the form $\sum a(n)w^n$ satisfying the ASD
congruence relations
$$a(np)-A(p)a(n)+B(p)a(n/p)=0 \mod p^{(k-1)(1+\ord_p n)}, \quad
\forall n\ge 1$$ where $A(p)$ and $B(p)$ are $p$-adic and in general
such bases are different for each $p$.
\smallskip

\subsection{}To examine the ASD congruence relations computationally,
Atkin uses ``$p$-adic" Hecke operators as. For any ``good" prime $p$, define the $p$-adic Hecke
operator $T_p=T^2+s(p)\cdot p^{k-1}$ for $ f=\sum a(n)w^n$  as
follows:
$$ T_p(f)= \sum_n (a(np)+s(p)p^{k-1}a(n/p) )  w^n.$$  Here, $s(p)$ is an algebraic number with absolute
value 1. In general, it is totally nontrivial to guess the right
$s(p)$ values since they are rarely Dirichlet characters.
 It is called a ``$p$-adic" Hecke operator as the $n$th
coefficient of $f|T_p$ will be determined up to modulo some suitable
power of $p$ depending on $n$. Unlike the classical Hecke operators
which can be diagonalized simultaneously, these operators are a
priori defined over various $p$-adic fields. It is clearly
exceptional that the space $S_k(\G)$ can be decomposed
simultaneously for almost all ``$p$-adic" Hecke operators.

\smallskip

\subsection{}In  \cite{sch85b}, Scholl has proven a collective version of Atkin and
Swinnerton-Dyer's observation. Basically, he showed that for any
form $f$ with rational coefficients  in a $d$-dimensional space of
even weight noncongruence cuspforms,  the coefficients of $f$
satisfy a $2d+1$ term congruence recursion. In \cite{L5}, the
authors have verified the necessary conditions under which Scholl's
results are valid for odd weight cases.

\subsection{}The available operators remain  the key players
when one attempts to reduce Scholl's result to 3-term ASD
congruence. In \cite{lly05, all05, long061} significant efforts have
been made to study the relations between these operators and
elements in the Galois groups. With the aid of these operators,
stronger results such as the existence of simultaneous (resp.
semi-simultaneous) eigen basis for  Atkin and Swinnerton-Dyer
congruences at almost all primes $p$ can be achieved cf.
\cite{lly05} (resp. \cite{all05,long061}).
\begin{example}[\cite{lly05}]
  There exists a basis $f_{\pm}=\sum a_{\pm}(n) w^n$ with coefficients in $ \Z[1/3](i)$ for $S_3(\G_3)$ and a pair of
  newforms $g_{\pm}=\sum b_{\pm} (n)q^n$ with coefficients in $\Z[i]$ such
  that for any $n\ge 1$ and any odd prime $p$ different from 3
  $$\frac{ a_{\pm}(np)-b_{\pm}(p) a_{\pm}(n)+\left ( \frac{ -3 }{p} \right )
  p^2 a_{\pm} (n/p)}{(np)^2}$$ is integral for some prime in $\Z[i]$
  above $p$.
\end{example}

In general, the symmetries inherited from $\G^N/\G$ is not
sufficient to derive 3-term Atkin and Swinnerton-Dyer congruences,
especially when the dimension $d$ is large \cite{alO07}. It is still
widely open whether the original 3-term Atkin and Swinnerton-Dyer
congruences always hold for any space of noncongruence cuspforms.


\newcommand{\etalchar}[1]{$^{#1}$}
\providecommand{\bysame}{\leavevmode\hbox
to3em{\hrulefill}\thinspace}
\providecommand{\MR}{\relax\ifhmode\unskip\space\fi MR }
\providecommand{\MRhref}[2]{%
  \href{http://www.ams.org/mathscinet-getitem?mr=#1}{#2}
} \providecommand{\href}[2]{#2}

\end{document}